\documentclass[10pt]{amsart}
\usepackage{amssymb}
\usepackage{epsfig}
\usepackage{url}
\usepackage{setspace}
\theoremstyle{plain}

\newtheorem{thm}{Theorem}[section]
\newtheorem{cor}[thm]{Corollary}

\newtheorem{rem}[thm]{Remark}
\newtheorem{ques}[thm]{Question}

\numberwithin{equation}{section}

\def\cal{\mathcal}
\def\bbb{\mathbb}
\def\op{\operatorname}
\renewcommand{\phi}{\varphi}

\newcommand{\Z}{\bbb{Z}}
\newcommand{\Q}{\bbb{Q}}

\newcommand{\C}{\bbb{C}}

\begin{document}
\title{Rational points on certain del Pezzo surfaces of degree one}
\author{Maciej Ulas}



\date{}

\begin{abstract}
Let $f(z)=z^5+az^3+bz^2+cz+d \in \Z[z]$ and let us consider a del
Pezzo surface of degree one given by the equation
$\cal{E}_{f}:\;x^2-y^3-f(z)=0$. In this note we prove that if the
set of rational points on the curve
$E_{a,\;b}:Y^2=X^3+135(2a-15)X-1350(5a+2b-26)$ is infinite, then
the set of rational points on the surface $\cal{E}_{f}$ is dense
in the Zariski topology.

\end{abstract}

\maketitle

\section{Introduction}\label{sec1}
A projective and geometrically irreducible surface $S$ is called a
del Pezzo surface if its anticanonical class $-K_{S}$ is ample.
The degree of del Pezzo surface is the self-intersection number of
its canonical class: $\op{deg}S=K_{S}^{2}$. As we know, the number
$\op{deg}S$ is positive integer and in fact we have inequality
$1\leq \op{deg}S \leq 9$. Geometrically, smooth del Pezzo surfaces
are obtained by blowing up $d\leq 8$ points in general position in
$\mathbb{P}^{3}$. The singular ones are blow-ups of
$\mathbb{P}^{2}$ in special configurations of points or in
infinitely near points. A detailed study of geometric and
arithmetic properties of del Pezzo surfaces can be found in
\cite{Manin}.

Many interesting arithmetic questions are connected with the class
of del Pezzo surfaces. As such surfaces are geometrically rational
(they are rational over the field $\C$), it is especially
interesting to look at problems concerning the question about
density (in Zariski topology) of $k$-rational points, where $k$ is
fixed finite extension of the field of rational numbers. From the
arithmetic and geometric points of view higher degree del Pezzo
surfaces are simpler than low degree ones. For example, it turns
out that in case $\op{deg}S=5$ the surface $S$ is always
$k$-rational. In case $\op{deg}S\leq 4$, no equivalent general
theorems are known.

In this note we are interested in the problem of existence of
rational points on del Pezzo surface of degree one given by the
equation
\begin{equation}\label{R1}
\cal{E}_{f}:\;x^2-y^3-f(z)=0,
\end{equation}
where $f(z)=z^5+az^3+bz^2+cz+d \in \Z[z]$. It is clear that after
a suitable change of coordinates each surface given by the
equation $px^2+qy^3+g(z)=0$, where $p,\;q\in
\Q\setminus\{0\},\;g\in \Q[z]$, $\op{deg}g=5$ or $\op{deg}g=6$ and
$g$ has a rational root can be transformed to the surface
$\cal{E}_{f}$ for a suitably selected $f\in\Z[z]$ of degree $5$.

Now, if we fix $z\in \Q$ such that $f(z)\neq 0$, then the curve
$\cal{E}_{f(z)}$ is an elliptic curve. So, in a natural way we can
consider $\cal{E}_{f}$ as on an elliptic surface where base curve
of the fibration is $z$-line. In our case $\cal{E}_{f}$ is an
isotrivial elliptic surface with $j$-invariant $j(\cal{E}_{f})=0$.
In recent years several papers have appeared concerning the
density of points in certain classes of isotrivial elliptic
surfaces. Interested readers can take a look at papers
\cite{KuwataWang},\;\cite{Rohrlich},\;\cite{Manduchi},\;\cite{BogomolovTschinkel},\;\cite{Munshi}.

In Section \ref{sec2} we prove Theorem \ref{thm1} which says that
if the set of rational points on the curve
\begin{equation*}
E_{a,\;b}:\;Y^2=X^3+135(2a-15)X-1350(5a+2b-26)
\end{equation*}
is infinite, then rational points on the surface $\cal{E}_{f}$ are
dense in Zariski topology. We also show that if $t$ is a
transcendental parameter, then under the same assumption on the
curve $E_{a,\;b}$ the set of $\Q[t]$-points on the surface
 \begin{equation*}
 \cal{E}_{f}:\;x^2-y^3-f(z)=t,
\end{equation*}
which is treated over the field of rational functions $\Q(t)$, is
infinite.

In Section \ref{sec3} we show as well that for every pair $a,\;b$
of nonzero rational numbers the set of rational points on the
surface
\begin{equation*}
 S:\;x^2+ay^5-z^6=b
\end{equation*}
(which is not a del Pezzo surface) is infinite.

\bigskip

\section{Rational points on the surface $\cal{E}_{f}$}\label{sec2}

Let $f(z)=z^5+az^3+bz^2+cz+d \in \Z[z]$ and let us consider the
surface
\begin{equation*}
\cal{E}_{f}:\;x^2-y^3-f(z)=0.
\end{equation*}

For a given $z\in\Q$ denote the curve $Y^2=X^3+f(z)$ by $E_{z}$.
Let us recall what the torsion part of the curve $E_{z}$ looks
like with a fixed $z\in\mathbb{Q}$  (\cite{Silverman}, page 323).
If $f(z)=1$, then $\op{Tors} E_{z}\cong\mathbb{Z}/6\mathbb{Z}$. If
$f(z)\neq 1$ and $f(z)$ is a square in $\Q$, then $\op{Tors}
E_{z}=\{\mathcal{O},\;(0,\;\sqrt{f(z)}),\;(0,\;-\sqrt{f(z)})\}$.
In the case that $f(z)=-432$ we have
$\op{Tors}E_{z}=\{\mathcal{O},\;(12,\;36),\;(12,\;-36)\}$. If
$f(z)\neq 1$ and $f(z)$ is a cube in $\Q$, then $\op{Tors}
E_{z}=\{\mathcal{O},\;(-\sqrt[3]{f(z)},\;0)\}$. In the remaining
cases we have $\op{Tors} E_{z}=\{\mathcal{O}\}$. The above have
two important consequences. If polynomial $f$ does not have
multiple roots, then there are only finitely many rational numbers
$z$ for which the curve $E_{z}$ has non-trivial torsion points.
Indeed, because every curve $f(z)=\lambda v^{i}$ where
$i=2,\;3,\;6$ and $\lambda\in\{1,-432\}$ is of genus $\geq 2$, so
our observation is an immediate consequence of Faltings theorem
\cite{Faltings}. Furthermore, if there is a rational base change
$z=\psi(t)$ such that on the surface $\cal{E}_{f\circ \psi}$ we
have the section $\sigma=(x,\;y)\in \Q(t)\times \Q(t)$ with
$xy\neq 0$, then $\sigma$ is a non-torsion section. It is clear
that we have the same conclusion in case that $f$ has rational
multiple root and $f(\psi(t))$ is nonconstant six-th power free
element of $\Q(t)$.

Now we are ready to prove the following

\begin{thm}\label{thm1}
Let $f(z)=z^5+az^3+bz^2+cz+d \in \Z[z]$ and consider the surface
given by the equation $\cal{E}_{f}:\;x^2-y^3-f(z)=0$. Then
\begin{enumerate}
\item if $f$ has multiple roots over $\C$, then there exists a
rational change of base $z=\psi(t)$ such that there is a
non-torsion section on elliptic surface $\cal{E}_{f\circ\psi}$. In
this case the set of rational points on $\cal{E}_{f}$ is Zariski
dense. \item if $f$ does not have multiple roots and the set of
rational points on the curve
$E_{a,\;b}:\;Y^2=X^3+135(2a-15)X-1350(5a+2b-26)$ is infinite, then
the set of rational points on the surface $\cal{E}_{f}$ is Zariski
dense.
\end{enumerate}
\end{thm}
\begin{proof}
(1) We have to consider two cases when $f$ has a rational multiple
root, and when $f$ has irrational multiple root.

If polynomial $f$ has a rational multiple root, then without loss
of generality we can assume that $f(z)=z^2(z^3+a'z^2+b'z+c')$ for
certain integers $a',\;b',\;c'$. For the proof of our theorem it
will be convenient to work with the surface $\cal{E}'_{f}$ given
by the equation
\begin{equation*}
\cal{E}_{f}':\;F'(X,Y,Z):=X^2-ZY^3-(Z^3+a'Z^2+b'Z+c')=0,
\end{equation*}
which is birationally equivalent to $\cal{E}_{f}$ by the mapping
$(x,\;y,\;z)=(ZX,\;ZY,\;Z)$ with the inverse
$(X,\;Y,\;Z)=(x/z,\;y/z,\;z)$. In order to find the rational curve
on $\cal{E}_{f}'$ we use method of indetermined coefficients. So,
we suppose that $X=Z^2+pZ+q,\;Y=Z+t$, where $t$ is transcendental
parameter and we are looking for $p,\;q \in \Q(t)$. For $X,\;Y$
defined in this way we obtain
\begin{equation*}
F'(X,Y,Z)=f_{0}+f_{1}Z+f_{2}Z^2+f_{3}Z^3,
\end{equation*}
where
\begin{equation*}
\begin{array}{ll}
  f_{0}=-c' + q^2, & f_{1}=-b'+2pq-t^3, \\
  f_{2}=-a'+ p^2+2q-3t^2, & f_{3}=-1 + 2 p - 3 t. \\
\end{array}
\end{equation*}

Let us note that the system of equations $f_{2}=f_{3}=0$ has
exactly one solution in $\Q(t)$ given by
\begin{equation}\label{R2}
p=\frac{1 + 3 t}{2},\quad q=\frac{-1+4a'-6t+3t^2}{8}.
\end{equation}
For $p,\;q$ defined in this way we can now solve the equation
$f_{0}+f_{1}Z=0$ according to $Z$. Then we obtain
\begin{equation*}
Z=-\frac{9t^4-36t^3+6(4a'+5)t^2-12(4a'-1)t+16a'^2-8a'-64c'+1}{8(t^3-15t^2+3(4a'-3)t+4a'-8b'-1)}=:\psi(t).
\end{equation*}
Therefore, we see that on the surface $\cal{E}_{f\circ \psi}$ we
have the section $\sigma
=(x,\;y)=(\psi(t)(\psi(t)^2+p\psi(t)+q),\;\psi(t)+t)$, where
$p,\;q$ are given by (\ref{R2}). With the help of computer we find
that $f(\psi(t))$ is nonconstant and six-th power free element of
$\Q(t)$, thus from the remark presented at the beginning of the
section it results that the section $\sigma$ is a non-torsion
section, so after the base change we get infinitely many sections
(corresponding to $m\sigma$), each with infinitely many rational
points. That means that each section is included in the Zariski
closure, say $\cal{R}$, of the set of rational points. Because
this closure consist of only finitely many components, it has
dimension two, and as the surface is irreducible, $\cal{R}$ is the
whole surface. Thus the set of rational points on $\cal{E}_{f\circ
\psi}$ is dense. The map from the $\cal{E}_{f\circ \psi}$ to the
$\cal{E}_{f}$ is dominant, so the set of rational points on the
surface $\cal{E}_{f}$ is also dense.

\bigskip

In the case that $f$ has irrational multiple root, without loss of
generality we can assume that $f(z)=(z^2+a')^2(z+b')$ for certain
integers $a',\;b'$ and $a'\neq 0$. After change of variables
$(x,\;y,\;z)=((Z^2+a')X,\;(Z^2+a')Y,\;Z)$ with the inverse
$(X,\;Y,\;Z)=(x/(z^2+a'),\;y/(z^2+a'),\;z)$ we reduce our problem
to the examination of the surface $\cal{E}_{f}''$ given by the
equation
\begin{equation*}
\cal{E}_{f}'':\;X^2-(Z^2+a')Y^3-(Z+b')=0.
\end{equation*}
Now, we look at the fibration of $\cal{E}_{f}''$ over the $Y$
line. We consider its generic fiber, which is a curve over the
field $\Q(Y)$, and viewing this curve over the extension field
$\Q(u)$ of $\Q(Y)$ given by $Y=u^2$. The equation of this curve is
of the form
\begin{equation*}
C:\;X^2=u^6Z^2+Z+a'u^6+b'.
\end{equation*}
Let us note that $C$ is a curve of genus zero with the $\Q(u)-$
rational point $P=(u^3,0,1)$ (it is a point at infinity). So, this
is a rational curve over $\Q(u)$. Putting $X=pu^3+t,\;Z=p$ and
solving the obtained equation according to $p$ we get the
parametrization of our curve
\begin{equation}\label{R3}
X(t,u)=\frac{u^3t^2-u+t^3(a't^6+b')}{2u^3t-1},\quad
Z(t,u)=\frac{-t^2+a'u^6+b'}{2u^3t-1}.
\end{equation}
Our reasoning show that the surface  obtained by the base change
$Y=u^2$ is in fact rational and thus that the surface
$\cal{E}_{f}''$ is unirational. This immediatly implies that
rational points on $\cal{E}_{f}$ are dense in the Zariski
topology.

\bigskip

(2) Let us put $F(x,y,z)=x^2-y^3-f(z)$. As in the previous proof
we use the method of indetermined coefficients. We look for
elements $p,\;q,\;r,\;T\in \Q(s,u)$ such that the curve given by
$x=T^3+pT^2+qT+r,\;y=T^2+sT+u,\;z=T$ lies on the surface . For
$x,\;y,\;z$ defined in this way we have
\begin{equation*}
F(x,y,z)=f_{0}+f_{1}T+f_{2}T^2+f_{3}T^3+f_{4}T^4+f_{5}T^5,
\end{equation*}
where
\begin{equation*}
\begin{array}{ll}
  f_{0}=-d+r^2-u^3, & f_{1}=-c+2qr-3su^2,  \\
  f_{2}=-b+ q^2 + 2pr - 3s^2u - 3u^2, & f_{3}=-a+ 2pq + 2r - s^3 - 6su, \\
  f_{4}=p^2 + 2q - 3s^2-3u, & f_{5}=-1 + 2 p - 3 s. \\
\end{array}
\end{equation*}

Let us notice that the system of equations $f_{3}=f_{4}=f_{5}=0$
has exactly one solution in $\Q(s,u)$ given by
\begin{equation*}
p=\frac{1+3s}{2},\; q=\frac{-1-6s+3s^2+12u}{8},\; r=\frac{1
+8a+9s+15s^2-s^3-12u+12su}{16}.
\end{equation*}
Putting these values into the equation $f_{2}=0$ we get
\begin{equation*}
48u^2-24(-3-10s+s^2)u-(5+32a-64b+60s+96as+198s^2+140s^3-3s^4)=0.
\end{equation*}
Solving now this equation with respect to $u$ we obtain
\begin{equation*}
u_{1,\;2}=\frac{-9 - 30s + 3s^2 \pm 4\sqrt{15s^3+
90s^2+9(2a+5)s+6(a-2b+1)}}{12}.
\end{equation*}

Now let us consider the curve $C_{a,\;b}$ given by the equation
\begin{equation*}
C_{a,\;b}:\;v^2=15s^3+ 90s^2+9(2a+5)s+6(a-2b+1).
\end{equation*}
After the affine change of variables
\begin{equation*}
(s,\;v)=\Big(\frac{X-30}{15},\;\frac{Y}{15}\Big)\quad \mbox{with
the inverse}\quad (X,\;Y)=(15(s+2),\;15v),
\end{equation*}
we transform the curve $C_{a,\;b}$ into the curve $E_{a,\;b}$
given by the equation
\begin{equation*}
E_{a,\;b}:\;Y^2=X^3+135(2a-15)X-1350(5a+2b-26).
\end{equation*}

It is an immediate consequence of our reasoning that in case that
 infinitely many rational points lie on $E_{a,\;b}$, then
all but finitely many points give us new value $z=T$ (the solution
of the equation $f_{0}+f_{1}T=0$) for which the elliptic curve
$\cal{E}_{f(z)}:x^2=y^3+f(z)$ has a positive rank. We know that
each fiber that contains infinitely many rational points is
contained in the closure of the set of rational points, say
$\cal{R}$. We have to show that for all but finitely many $z$ the
curve $\cal{E}_{f(z)}$ is of positive rank, so the set $\cal{R}$
contains infinitely many curves and is therefore 2-dimensional. As
the surface is irreducible, the set of rational points is dense.


\end{proof}

\begin{rem}\label{rem1}{\rm In the first part of Theorem \ref{thm1}
we have shown that if the polynomial $f$ has a multiple root, then
there is a rational change of base $z=\psi(t)$ such that there
exists a non-torsion section on the elliptic surface
$\cal{E}_{f\circ\psi}$. A natural question arises whether it is
possible to construct the polynomial $f$ without multiple roots
which would give a surface $\cal{E}_{f}$ with a similar property.
It turns out that it is possible to construct demanded
polynomials. This is tightly connected with the question for which
rational numbers $a,\;b$ the curve $E_{a,\;b}$ is singular.

 The curve $E_{a,\;b}:\;Y^2=X^3+135(2a-15)X-1350(5a+2b-26)=:g(x)$
 is singular if for the discriminant of the polynomial $g$ we have
\begin{equation*}
\Delta(a,b)=-16(4(135(2a-15))^3+27(1350(5a+2b-26))^2)=0.
\end{equation*}
 It is easy to see that it holds if and only if
 \begin{equation*}
45(2a-15)=-t^2,\quad \quad 675(5a+2b-26)=t^3,
 \end{equation*}
 for a certain rational number $t$. Hence we obtain
\begin{equation}\label{R4}
a=\frac{-t^2+630}{90},\quad \quad b=\frac{2t^3+75t^2-15525}{2700}.
\end{equation}
For $a,\;b$ defined in this way, the curve $E_{a,\;b}$ is reduced
to the form $E:\;Y^2=(X+2t)(X-t)^2$ and the set of its rational
points can be parametrized in the following way
 \begin{equation*}
X=U^2-2t,\quad Y=U(U^2-3t).
 \end{equation*}
 Using the above equalities we obtain an explicit form
 of the rational curve  $L$ on the surface $\cal{E}_{f}$ in case when
 $a,\;b$ are given by (\ref{R4}). We will not, however, present these
 equations here since the rational functions defining the curve $L$ are of
 very high degrees.}
 \end{rem}

Part (1) of our theorem and the above observation suggest the
following

\begin{ques}\label{ques1}
Let $f(z)=z^5+at^3+bz^2+cz+d$ and let us consider the surface
$\cal{E}_{f}:\;x^2=y^3+f(z)$. What are the conditions guaranteeing
the existence of a rational base change $z=\psi(t)$ such that
there is a non-torsion section on the surface
$\cal{E}_{f\circ\psi}$?
\end{ques}

Let us note several interesting corollaries of Theorem \ref{thm1}.

\begin{cor}\label{cor1}
Let $f(z)=z^5+az^3+bz^2+cz+d\in\Z[t]$ and let us assume that the
set of rational points on the curve
$E_{a,\;b}:\;Y^2=X^3+135(2a-15)X-1350(5a+2b-26)$ is infinite.
Then, the diophantine equation $x^2-y^3-f(z)=t$ has infinitely
many solutions in the ring of polynomials $\Q[t]$.
\end{cor}
\begin{proof}
This observation is an immediate consequence of the reasoning
conducted in the proof of the second part of \ref{thm1}. Indeed, by
the assumption concerning the curve $E_{a,\;b}$ we have shown that
for fixed rational point on $E_{a,\;b}$ the coefficients of the
polynomials $x(T)=T^3+pT^2+qT+r,\;y(T)=T^2+sT+u,\;z(T)=T$ are
rational and the polynomial $F(x(T),\;y(T),\;z(T))$, where
$F(x,y,z)=x^2-y^3-f(z)$, is of degree one. Because the set of
rational points on the curve $E_{a,\;b}$ is infinite, then for all
but finitely many points from $E_{a,\;b}(\Q)$ we have $f_{1}\neq 0$.
Now, solving the equation $f_{0}+f_{1}T=t$, we get that
$T=(t-t_{0})/f_{1}$. Finally we obtain the identity
\begin{equation*}
x\Big(\frac{t-t_{0}}{f_{1}}\Big)^2-y\Big(\frac{t-t_{0}}{f_{1}}\Big)^3-f\Big(\frac{t-t_{0}}{f_{1}}\Big)=t.
\end{equation*}
This concludes the proof of Corollary \ref{cor1}.
\end{proof}

\begin{cor}\label{cor2}
Let $c,\;d\in \Z$ and $f(z)=z^5+cz+d$. Then on the surface
$\cal{E}_{f}$ the set of rational points is dense in the Zariski
topology.
\end{cor}
\begin{proof}
In case $a=b=0$ the equation of the curve $E_{a,\;b}$ from the
proof of Theorem \ref{thm1} takes the form
\begin{equation*}
E_{0,\;0}:\;Y^2=X^3-2025X+35100.
\end{equation*}
With the assistance of {\sc APECS} program \cite{Connell} we have
found that $\op{Tors}(E_{0,\;0})$ is trivial and that the rank of
the curve $E_{0,\;0}$ is two. Independent points of infinite
order, $P_{1}=(15,90)$ and $P_{2}=(25,10)$, generate the set of
rational points on $E_{0,\;0}$. It means that the set of rational
points on $E_{0,\;0}$ is infinite. This ends the proof of our
corollary.
\end{proof}

\begin{cor}\label{cor3}
Let $a,\;b,\;c,\;d\in \Z\setminus\{0\}$ and consider the surface
given by the equation
\begin{equation*}
\cal{E}:\;ax^2+by^3+cz^5=d.
\end{equation*}
Then, the set of rational points on $\cal{E}$ is dense in the
Zariski topology.
\end{cor}
\begin{proof} After the change of variables
$x=X/a^8b^{10}c^{12},\; y=-Y/a^5b^7c^8,\; z=-Z/a^3b^4c^5$, the
equation of the surface $\cal{E}$ takes the form
\begin{equation*}
X^2-Y^3-Z^5=a^{15}b^{20}c^{24}d.
\end{equation*}
From Corollary \ref{cor2} we get the statement of our theorem.

For instance, if we take point $P_{1}=(15,\;90)$, which is one of
generators of the set $E_{0,\;0}(\Q)$ and after we perform all
necessary calculations presented in the proof of Theorem
\ref{thm1} we obtain the identity
\begin{align*}
&a\Big(\frac{25875323c^{18}+720748a^{15}b^{20}c^{12}d+8336a^{30}b^{40}c^6d^2+64a^{45}b^{60}d^3}{1560896a^8b^{10}c^{15}}\Big)^2+\\
&b\Big(-\frac{87709c^{12}+1544a^{15}b^{20}c^6d+16a^{30}b^{40}d^2}{13456a^5b^7c^{10}}\Big)^3+c\Big(\frac{135c^6+4a^{15}b^{20}d}{116a^3b^4c^5}\Big)^5=d.
\end{align*}
Note that, in fact, we proved something more. Namely, for every
natural number $d$  there is a $S$-integer point on the surface
$\cal{E}$, where $S=\{p\in \mathbb{P}: \;p|58abc\}$.
\end{proof}

\bigskip

\section{Rational points on the surface
$x^2+ay^5-z^6=b$}\label{sec3}

In this section we will show that the method employed to prove
Theorem \ref{thm1} can be used in another situations. Interested
readers can also take a look at paper \cite{Ulas}. We will prove
the following

\begin{thm}\label{thm2}
For every pair $a,\;b$ of nonzero rational numbers the set of
rational points on the surface $\cal{S}:\;x^2+ay^5-z^6=b$ is
infinite.
\end{thm}
\begin{proof}
Our reasoning will be similar to that presented in the proof of
Theorem \ref{thm1}. Let us denote $F(x,y,z)=x^2+ay^5-z^6$ and let
us put $x=T^3+pT^2+qT+r,\;y=uT+v,\;z=T$. For $x,\;y,\;z$ defined
in this way we have
\begin{equation*}
F(x,y,z)=f_{0}+f_{1}T+f_{2}T^2+f_{3}T^3+f_{4}T^4+f_{5}T^5,
\end{equation*}
where
\begin{equation*}
\begin{array}{lll}
  f_{0}=r^2+av^5, & f_{1}=2qr+5auv^4,  &  f_{2}=q^2+2pr+10au^2v^3,\\
  f_{3}=2pq+2r+10au^3v^2,  &  f_{4}=p^2+2q+5au^4v, & f_{5}=2p+au^5.  \\
\end{array}
\end{equation*}

The system of equations $f_{2}=f_{3}=f_{4}=f_{5}=0$ has exactly
three solutions (with respect to $p,\;q,\;r,\;v$). One of them is
defined over $\Q(u)$ (the other two are defined over
$\Q(\sqrt{3})(u)$), namely
\begin{equation}\label{R5}
p=-\frac{au^5}{2},\quad q=\frac{3a^2u^{10}}{16},\quad
r=\frac{a^3u^{15}}{64},\quad v=-\frac{au^6}{8}.
\end{equation}
Now putting the calculated values to the expressions defining
$x,\;y,\;z$, we see that the polynomial $F(x,y,z)$ is of degree
one. Solving the equation
\begin{equation*}
F(x(T),y(T),z(T))=b
\end{equation*} with the respect to $T$ and
performing all necessary calculations we obtain the identity
\begin{align*}
&\Big(\frac{118441a^{18}u^{90}+2^{15}11863a^{12}bu^{60}-2^{30}137a^{6}b^{2}u^{30}+2^{45}b^{3}}{2^{9}29^{3}a^{15}u^{75}}\Big)^2+\\
&a\Big(-\frac{9a^{6}u^{30}-2^{13}b}{58a^{5}u^{24}}\Big)^5-\Big(\frac{7a^{6}u^{30}-2^{15}b}{232a^{5}u^{25}}\Big)^6=b.
\end{align*}
This concludes the proof of Theorem \ref{cor2}.
\end{proof}

Using exactly the same reasoning as in the proof of the Theorem
\ref{thm2} we get

\begin{cor}\label{cor4}
Let $a,\;b,\;c,\;d\in \Z,\;a\neq 0$ and consider the surface given
by the equation $\cal{S}:x^2+ay^5+by-(z^6+cz)=d$. Then the set of
rational points on $\cal{S}$ is infinite.
\end{cor}

\bigskip

 \noindent Jagiellonian University, Institute of Mathematics, Reymonta 4, 30 - 059 Krak\'{o}w, Poland

 \noindent E-mail:\;{\tt Maciej.Ulas@im.uj.edu.pl}

 \end{document}